\documentclass[]{interact}

\usepackage{epstopdf}
\usepackage[caption=false]{subfig}

\usepackage[longnamesfirst,sort]{natbib}
\bibpunct[, ]{(}{)}{;}{a}{,}{,}

\theoremstyle{plain}

\theoremstyle{definition}

\theoremstyle{remark}

\begin{document}


\title{Creating and experiencing Flipped Learning in Multivariable Calculus for Engineering}

\author{
\name{Hugo Caerols-Palma\textsuperscript{a}\thanks{CONTACT Hugo Caerols-Palma Email: hugo.caerols@uai.cl} and Katia Vogt-Geisse\textsuperscript{a}\thanks{CONTACT Katia Vogt-Geisse Email: katia.vogt@uai.cl}}
\affil{\textsuperscript{a} Facultad de Ingenier\'ia y Ciencias, Universidad Adolfo Ib\'a\~nez, Diagonal Las Torres 2640, Pe\~nalol\'en, Santiago, Chile.}
}

%

\maketitle

\begin{abstract}

This article discusses the process of creating, implementing and experiencing Flipped Learning in a Multivariable Calculus course for second year engineering students. We describe the construction of the teaching material, consisting of short videos for pre-class preparation and aligned worksheets for in-class dynamics, and the activities that were conducted. We discuss difficulties and key aspects to be considered while creating this material and during implementation of Flipped Learning. We present how students reacted to pre-class preparation and how in-class dynamics developed during implementation.  
We show results on students performance and perception when enrolling in a flipped classroom section. We present comparative results on students performance of a section taught with Flipped Learning vs a parallel section thought in the traditional expository way. We could conclude that flipped courses show similar results in passing percentage than traditionally taught courses, that student's perceptions are generally mixed, and we perceived that students repeating the course preferably do not choose flipped classes.  Finally, we discussed the methodological  evolution of this course converging to a mixed methodology throughout a four year period, observing that the instructors evaluation decreases in classes that were flipped. Mixed methodologies on the other hand, increased the learning experience of students resulting in an   increased  instructors evaluation score and  higher students enrollment in the course. 

\end{abstract}


\section{Introduction}

In this paper we describe our experience in the preparation and implementation of an active learning methodology applied to the course \textit{Multivariable Calculus} for second year engineering students. In particular we describe:  the difficulties that appeared while creating the material and during implementation; the results about student's performance and perception upon implementation; the methodological evolution in the course in a four year period, and the evolution of indicators that reflect student's preferences for this active methodology during the same four years.

Active learning methodologies were introduced in the $90$'s by Eric Mazur, a physics Harvard professor, who first changed the way to transmit knowledge to university students. He based his teaching in active learning and peer instruction, incorporating among others the use of clickers and peer discussion in introductory physics courses at Harvard University. He observed increased student mastery in conceptual reasoning and quantitative problem solving in physics \citep{Mazur10,EM}. A key factor that makes active independent learning possible is technology. For instance, the use of video lectures to transmit knowledge has been widely used, for example by Salman Khan, an MIT alum, when he founded Khan Academy \citep{Kahn} in 2006. This online platform provides an open-source video library with thousands of videos on several topics, for instance college Calculus \citep{BV, Khan}. 



An active learning methodology that has been widely used is Flipped Learning (or sometimes referred to as Flipped Classroom), which is defined in a very simple way by \citet{DefFlipped1} as: \textit{``Inverting the classroom means that events that have traditionally taken place inside the classroom now take place outside the classroom and vice versa"}. We adopted the definition stated by the Flipped
Learning Network \citep{DefFlipped2}, since it describes in a more precise way what we understand by Flipped Learning, which is: \textit{``Flipped learning is a pedagogical
approach in which direct instruction
moves from the group learning
space to the individual learning
space, and the resulting group space
is transformed into a dynamic,
interactive learning environment
where the educator guides students
as they apply concepts and engage
creatively in the subject matter.''}
Additionally, we specify that direct instruction outside of class (pre-class preparation) is provided through video material, as is stated in the definition given by \citet{BV}.

Flipped Learning has firstly become popular for secondary education at a small Colorado highschool in the United States in 2007 \citep{Bergmann2012, Flaherty2015} and since then has  been applied also at university level courses in different ways. For instance, \citet{Wasserman(2015)}  studied the effect of flipping an undergraduate Calculus III class on student performance and perceptions. They found that, compared to a traditional expository class, student performance did not show any improvement in certain aspects and in others became at most slightly better, and that student perceptions are mixed. \citet{Sun2018} analyzed how students' self-regulatory factors such as self-efficacy in learning mathematics and help seeking strategies were related to academic achievement in a  flipped Calculus I and a Calculus II undergraduate course. Their results show that these factors significantly impact students' achievement in the flipped classroom setting. \citet{Toto2009} explored students' perception on how to improve the flipped classroom learning strategy to enhance learning in engineering education. Their results give insight into improvement of video material and in-class activities, for instance that videos should include more detail or that in-class activities should be designed in such a way as to keep students busy all the time.  \citet{Hao2016} investigated students' learning readiness with respect to Flipped Learning in two undergraduate courses of an education department. One of their results show that about $60\%$ of the students recognize the advantage of Flipped Learning, but less than $50\%$ stated that they would like to take a flipped class again, and that it seems that a significant number of students are not ready to take responsibility for their own learning. Some other examples of the use of Flipped Learning are listed in the reviews by \citet{Flaherty2015} and \citet{BV}.


During the years 2016 and 2017 we decided to explore incorporating Flipped Learning into the course of Multivariable Calculus at Adolfo Ib\'a\~nez University in Chile for second year undergraduate engineering students. Our flipped methodology was adopted from the definition stated above. I.e. it consisted in engaging students to learn independently trough short videos outside of the classroom as preparation for in-class activities, and time in class was used to work towards reaching the desired learning outcomes and proficiency, while actively involving students in knowledge construction. In class, students were coached by the professor while they worked actively on solving exercises, as opposed to being passive receptors of content like usually happens in a traditional expository class setting. In a flipped class, students apply, analyze and synthesize while guided by the professor. The audiovisual material was created by the professor who taught the class. 

The process was challenging in various ways: First, the design, elaboration and edition of the videos required us to gain specific technical knowledge; second, we had to decide on a student-centered learning theory for the in-class activities, which is a crucial part of implementing Flipped Learning \citep{BV}; third, we encountered difficulties during implementation that we needed to address; and fourth, we had to analyze in which way a methodological change of this type could endure in time and how we could give proper future use to the learning material we created.




 One of the reasons that led us to try this new  methodological approach is that the current trends in the curriculum design is pointing towards a curriculum based on measuring learning outcomes and abilities. In particular our  institution is ABET (Accreditation Board for Engineering and Technology \citep{ABET}) accredited for industrial engineering and hence must follow processes that guarantee that our graduates meet specific learning outcomes. These outcomes are naturally assessed and taught with the Flipped Learning methodology, since many ABET criterion are related to independent work, effective communication and learning-by-doing. Those aspects are included when we incorporate students actively in the learning process.
In that sense, we aimed to analyze how students react when basic mathematics courses aim to  contribute to a more integral education through Flipped Learning.

 


After initially obtaining support from specialists in education from our institution, we did not count on specialist educational developers or a technological team to assist continuously with the construction of the material. So we chose to flip the course Multivariable Calculus because of the affinity and experience of the professors teaching it and our interest in investing personal time and resources into this project. We were also aware of similar activities implemented at Purdue University in the United States for mathematics courses \citep{IC} because of connections of one of the authors, which gave us the idea to try and adapt new methodologies here in Chile. 

The novelty of our work consist of mainly two aspects: First,to the best of our knowledge, there is no study that describes in detail the experience of an instructor that went though the whole process, from thinking about implementing Flipped Learning, through creating the material and implementing the methodological change, up to obtaining students' perceptions and analyzing the evolution of the change in time. As \citet{Flaherty2015} mentioned, the papers in the literature usually study single course outcomes.
Our paper contributes to the literature since it describes the whole experience from the viewpoint of the instructor, providing information about difficulties and challenges that were encountered in the process. Second, cultural differences may determine many times the success or failure of different teaching methods, for instance \citet{Joanne2014} stated that flipped classroom was well-received among Asian students, but in other parts of the world often students' opinions are mixed \citep{BV, Flaherty2015}. There are little documented studies about the experience of implementing active methodologies, in particular Flipped Learning, in Chilean universities and how Chilean students may react to it. \citet{Rodriguez2017} described results on flipping some classes of a massive engineering course on \textit{Organizational Behavior} using a MOOC (massive open online course) platform. They describe among other aspects, that students showed more motivation when classes were flipped and more time should be left for discussion in class. \citet{Galido2016} discuss perceptions of secondary education students and teachers about Flipped Learning. To the best of our knowledge, there are no documented studies about implementing Flipped Learning in mathematics courses in Chile. 

This article is structured in the following way. In Section \ref{Material} we describe the process  of constructing the audiovisual material, the challenges encountered and the valuable experience we gained from it. In Section \ref{Pre-In-class} we explain the course work-flow, show how our students reacted to pre-class preparation using the video material, and give details about in-class activities commenting on how students behaved during that time. Section \ref{Results} presents  results on students' performance and perception, and the evolution of the methodology and students' preferences for Flipped Learning with time. Finally, in Section \ref{DisConcl} we discuss the results and give conclusions.


\section{Construction of the material }\label{Material}

Instead of using pre-existing video resources we decided to create our own video material for pre-class preparation. We thought this would be valued greatly by our students and would allow a closer teacher-student relationship and bigger student's compromise to watch them. This way we could also integrate in a better way pre-class activities with face-to-face class time according to what we considered key aspects for our students to learn. 
During the period of 2 semesters, we created 104 videos, a total of 21.43 continuous video hours. The duration of each crucially depended on how deep the ideas where treated in the video, even though we tried to keep the length at about 10 minutes.  Videos with only examples used to be shorter than others that show proofs. Table \ref{LengthVideos} classifies the videos according to their length. One can observe that the length of more than half of the videos is at most 12 minutes. 


\begin{table}[h!]
\centering
\tbl{Distribution of number of videos according to their duration in minutes.\\}
{\begin{tabular}{ccc} \midrule
 Duration $(t)$ in minutes & Number of videos of length $t$ & \% of videos of length $t$ \\ \midrule	
$ t < 8$ & 13 & 12.5  \\
$8\leq t<12$ & 47& 45.19 \\
$12\leq t <16$ & 21& 20.19 \\
$16\leq t <20$ &21& 20.19  \\
$ t > 20$& 2 & 1.92  \\ \bottomrule
\end{tabular}}
\label{LengthVideos}
\end{table}



It took us a long time to find the best way to record a video for the purpose of flipping our class. We performed a video-material review to get familiar with the type of video formats available before recording our own. We reduced our search to mainly two possibilities for video structures:  the first shows a lecturer explaining a topic on a blackboard (recorded while in class or without an audience) and the second shows only a screen where text appears while the voice of the lecturer explains the topic. There are plenty of videos on Youtube \citep{Youtube} and other platforms that serve as examples; for instance, the online course ``Understanding Einstein: The Special Theory of Relativity" from Stanford Online offered through Coursera   \citep{CourseraCourse, Co} uses the fist video structure, and Khan Academy \citep{Khan} videos use the second video structure. Variations of these two structures are generally left to the creativity of each lecturer, reflecting his/her teaching experience and methods with which they have succeeded in explaining a certain topic to their students. These methods vary from lecturer to lecturer. Some lecturers use graphic resources embedded in their videos, others animations and others just use their writing and plotting skills. 


After watching many examples of existing videos we realized key aspects that make a good pre-class video: first, good videos were carefully edited, the videos are immediately pleasant to watch and images and audio are high quality; second, it is noticeable that the video follows a structured script that develops gradually the difficulty of the the topic to be taught; third, examples are given that answer questions previously stated, guiding the student through the learning process. 



We chose to record our videos using the second video format, i.e. we do not appear physically on the video, only our writing is visible and our voice can be heard explaining the topic. We chose this approach because of mainly two reasons: first, the lack of resources to record more sophisticated videos and second we didn't need a special recording room or special recording conditions, we could record our videos at home or at the office using minimum resources. We used a special tablet, a drawing software and a screen capture software  that allows users to write simulating a white board resulting in good quality videos with not much editing knowledge   \citep{Wacom, Paint, Camtasia}. 

Writing on the tablet was relatively easy to learn, even though one can experience some initial difficulties such as achieving a good sized and nicely legible handwriting in order to use the screen size in optimal manner, and choosing font size and colors strategically to produce attractive videos. 
Creating our own videos that are of the quality and standard one desires and that transmit the material the way one had it in mind was also very time consuming. To find the best way to present a topic takes time, especially for professors with little experience in creating digital content. Whether to present a historical context, a graph, an animation and to embed these ideas into a script requires experience. Once a script is created it becomes easier to record the video. Initially a script can be as detailed as to include which colors the writing in the video should be showing and the layout of the ideas that will be presented. The script should then be followed while recording to obtain a quality video. Over time and after recording many videos, this process becomes easier and more fluent with experience, until sometimes managing to record without previously creating a  script. 

For some instructors, the learning curve to being able to create quality videos is long, while for others mastery is reached quite fast. 
To agree with the final product is a long process, for instance, getting used to talk and write simultaneously following a script may not be trivial as well as the editing process is not.
We experienced that not all professors feel comfortable creating digital teaching material. During a time period of two years, out of ten professors of the Mathematics Department of our Institution that started recording videos attracted by the novelty of it, only two kept recording enough videos to flip a class. The remaining professors did not continue because of several reasons, for instance, they didn't find the time to record and edit the videos, they didn't like the hardware or software that was used, they didn't like the quality of the videos they created, or they never felt comfortable recording their own videos. 


Despite of the difficulties we encountered, to record our own videos was a very valuable experience: 
First, the videos reveal weaknesses each professor may have while teaching, for instance, slow pace, monotone voice,  frequent pausing or a distracting use of fillers. To realize those weaknesses that may have passed unnoticed before, may help a professor improve his or her face-to-face teaching and interaction with students. One of the professors that participated in this project improved significantly his speaking, and became much more fluent in his traditional expository classes. Second, we realized the different ways that different professors explain the same topic. Even though mathematics is an exact science, when explaining a topic, the way the video develops and emphasizes different ideas varies greatly between professors. Watching videos created by colleagues enriches ones viewpoint on how to teach that specific topic treated in the video and interesting discussions arise on how to best present the material. Those discussions are very valuable for the education of instructors. In general, recording those videos allows the professors to be creative in how to present a topic. For instance, concepts can be introduced one after the other in a linear manner, or in a non-linear way, where the concepts are initially given and connected to each other afterwards (see two video examples in the supplementary material \ref{}). 
Third, while creating videos one needs to in parallel create material which will be used during class time for students to work based on the video content. This leads you to restructure the chronogram of the class, since by flipping the class  you rethink carefully and realize the real difficulty students may have with some of the topics and how long it may take to work those topics with them. As opposed to in a traditional expository class, where there are topics that are treated sometimes without much examples since the class is based solely on transmission of content, and additional time for students to  acquire deeper knowledge is not accounted for in class planning.  When active learning methodologies are used, the time used to cover subjects are often determined by the student's ability to reach the learning outcomes and not by the instructors ability to explain the content without student's active participation.

\section{Pre-class and in-class dynamics and experience}\label{Pre-In-class}

\begin{figure}[h!]
 \centering
\includegraphics[scale=0.33]{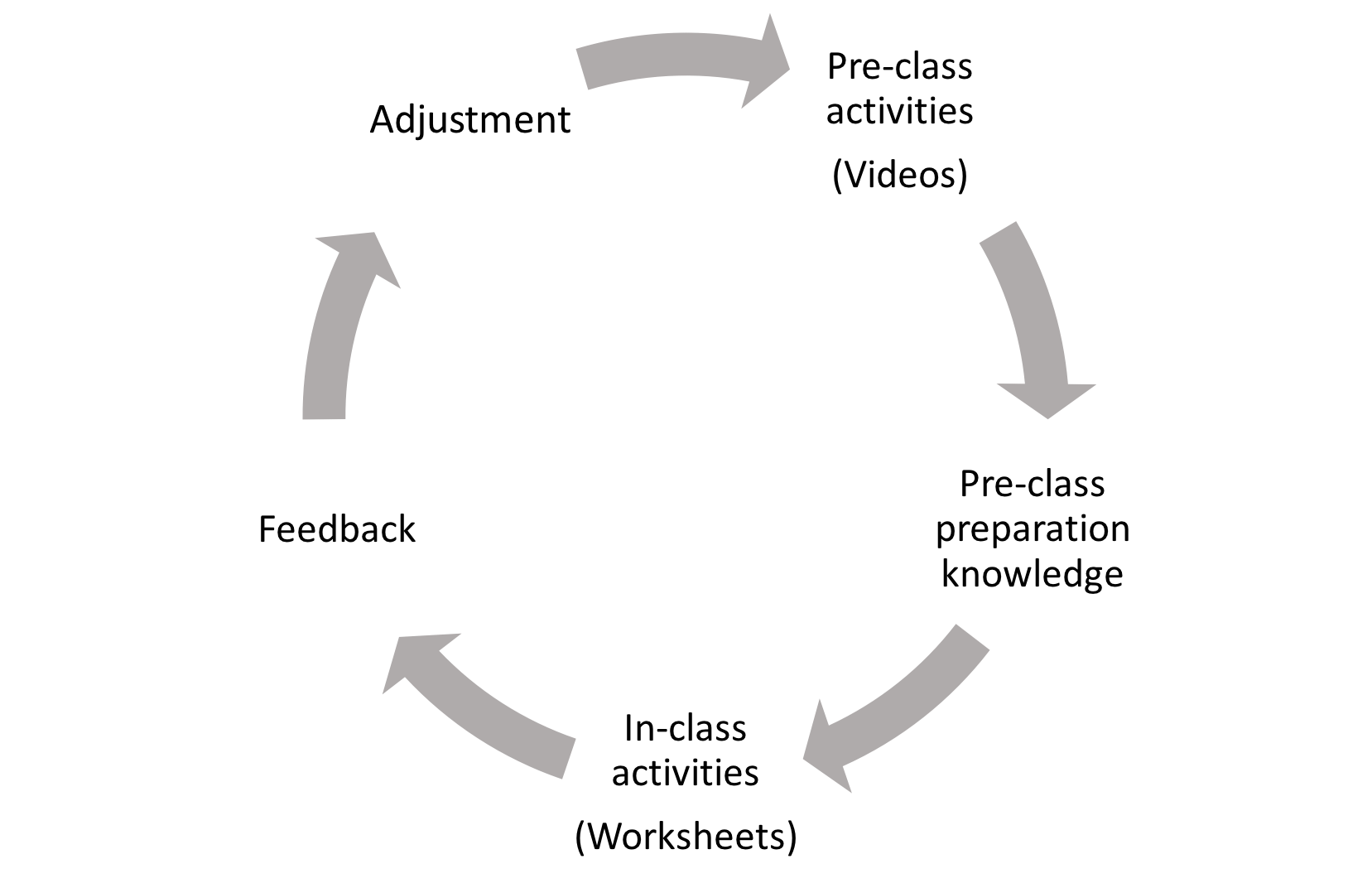}
\caption{The figure shows a diagram describing the working dynamics of the Flipped Learning methodology applied to the course of Multivariable Calculus throughout the semester.}
\label{Dynamics}
\end{figure}

In this section we describe how the flipped working dynamics was structured and experienced throughout the whole semester. Figure \ref{Dynamics} shows an explanatory diagram for the work flow.


 Students had to use time outside of the classroom for pre-class activities, which consisted in studying short videos. Those videos were posted on a Moodle \citep{Moodle} learning platform and were available to students anytime and anywhere through this platform a few days before class time.  While studying the videos, students acquired pre-class knowledge that prepared them to confront the in-class activities. How students confronted pre-class preparation and its effect on in-class dynamics is described in Subsection \ref{Pre-class}. Each class they were given a worksheet that was aligned with the videos they had watched, and they had to solve the worksheet collaboratively with their peers. While coached by the instructor during in-class time, students received feedback about their reasoning and the instructor received feedback about his or her students' progress, weaknesses and learning pace. Subsection \ref{In-Class} describes the in-class dynamics experience. Finally, the feedback obtained led to an adjustment period, where the professor adjusted the class pace and students should ideally have adjusted their pre-class learning time.

\subsection{Pre-class activities and preparation}\label{Pre-class}

\begin{figure}[htbp]
\centering
\includegraphics[scale=0.41]{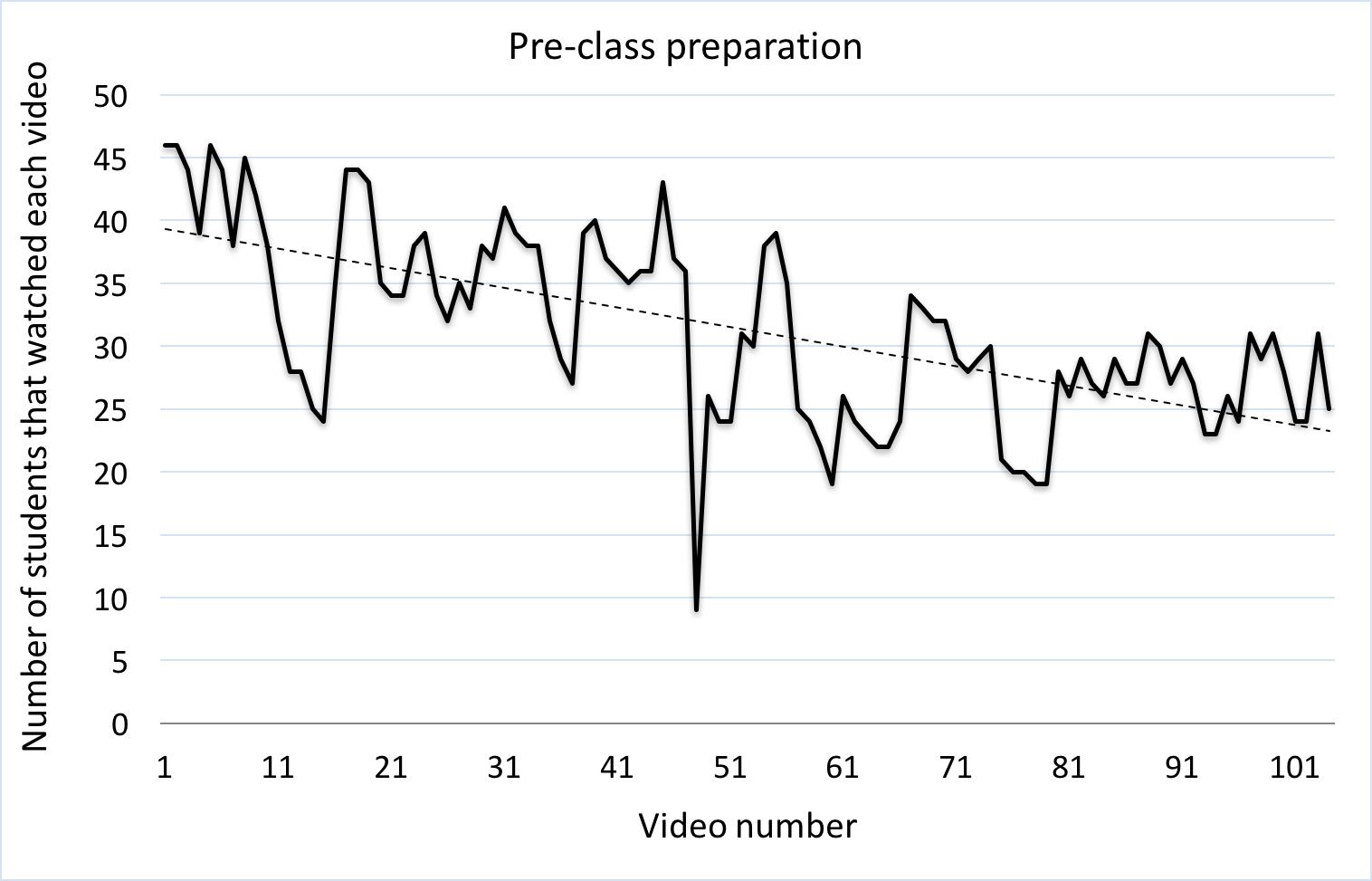}
\caption{The figure shows video number 1 to video number 105 that were posted before each class vs the number of students that watched each of those videos as pre-class preparation. The data corresponds to the implementation of Flipped Learining in Multivariable Calculus during the second semester of 2016 in a section with a total of 54 students enrolled.}
\label{PreClassPreparation}
\end{figure}

Our second year engineering students have a high course load. They enroll in six core courses each semester, each of them with a considerable high credit load. Throughout the years, not all students are efficient in organizing their time to meet all requirements for each course, and hence  we observed that not all students are able to watch the videos of a flipped class consistently as pre-class preparation. Figure \ref{PreClassPreparation} shows video number 1 to video number 105 vs the number of students that watched each of the videos as pre-class preparation, out of a total of 54 students enrolled in the course the second semester of 2016.
  The curve shows a decreasing trend along the semester. Initially, while videos 1 to 11 became available, students were eager to watch them and between 40 and 45 students did. As long as the semester progresses, several low peaks appear, which indicate low pre-class preparation. We know those coincided with weeks in which students were evaluated in other courses. By the end of the semester the curve tends to stabilize, being between 20 and 30 the number of students (around $50\%$ of the students) that watched the last videos (nr. 81 to 105).





\begin{figure}[htbp]
\centering
\includegraphics[scale=0.44]{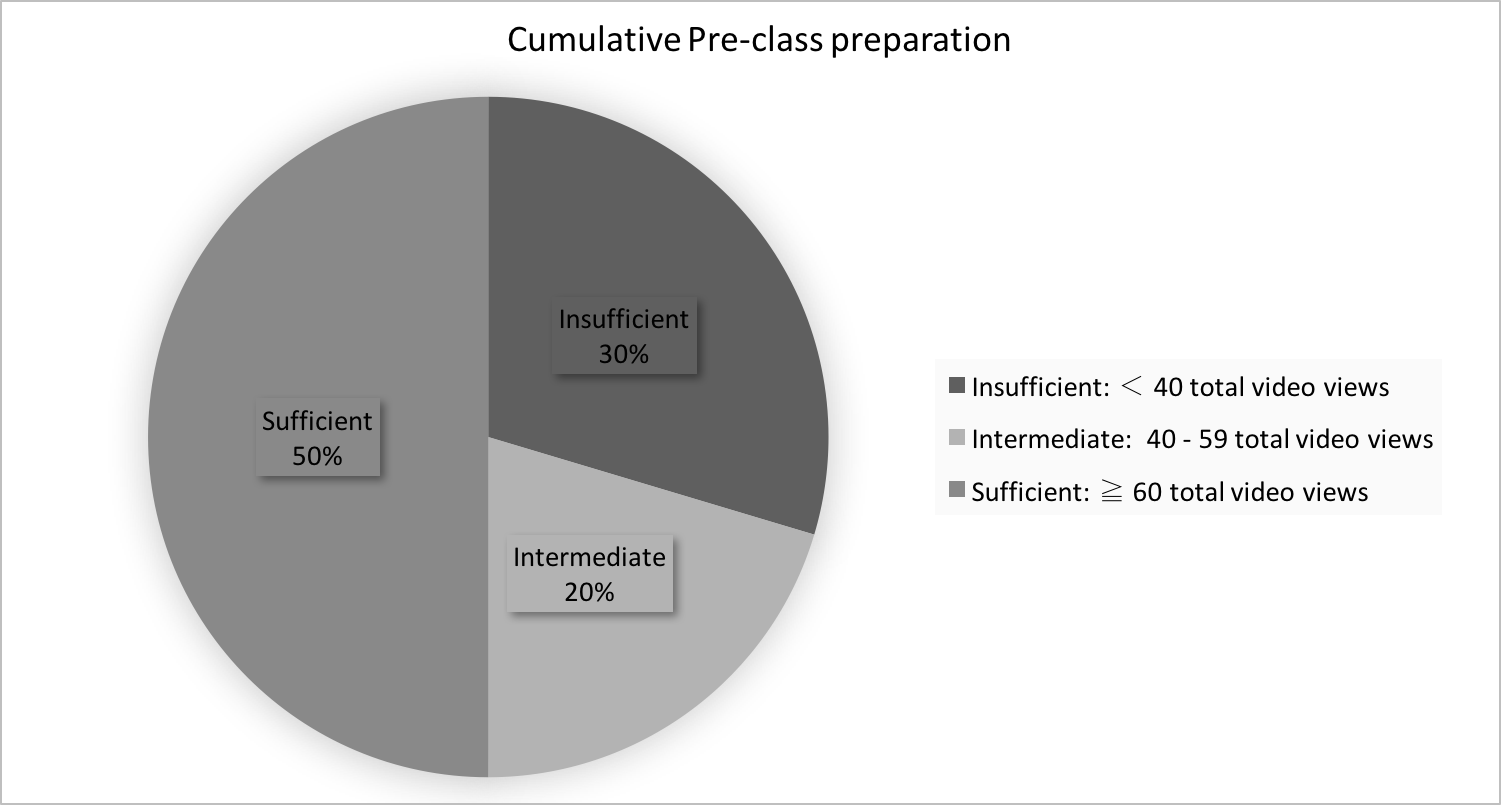}
\caption{The figure shows the cumulative pre-class preparation represented as the percentage of students that watched an insufficient (less than 40 videos), intermediate (between 40 and 59 videos), and a sufficient (at least 60 videos) number of videos as pre-class preparation throughout the whole semester.  The total number of videos that were available throughout the semester is 105. The data corresponds to the implementation of Flipped Learining in Multivariable Calculus during the second semester of 2016 in a section with a total of 54 students enrolled.}
\label{CumPreClassPreparation}
\end{figure}

\begin{figure}[h!]
\begin{center}
\includegraphics[scale=0.41]{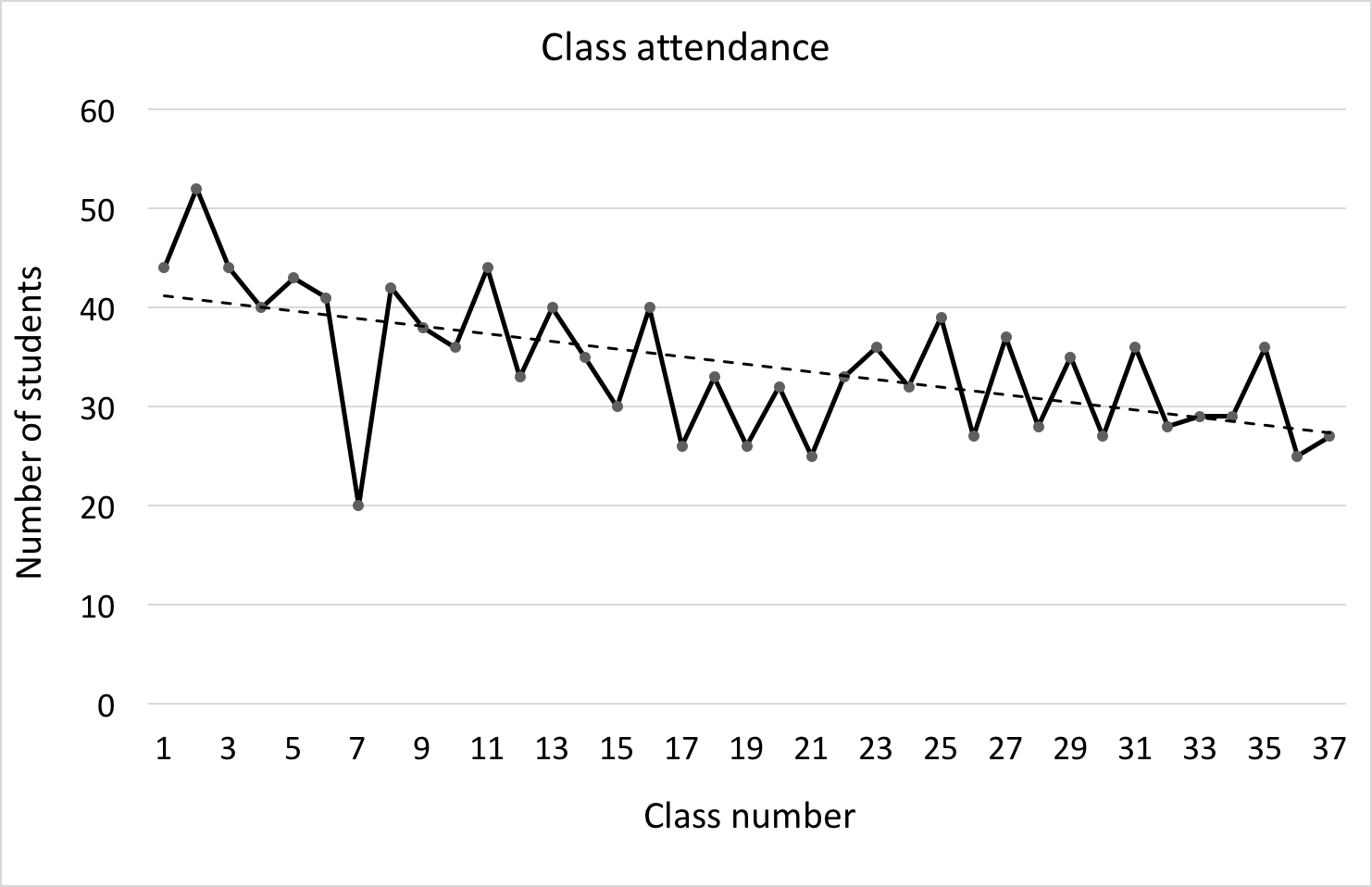}
\end{center}
\caption{The figure shows class numbers 1 to 41 vs number of students that attended each class.  The data corresponds to the implementation of Flipped Learining in Multivariable Calculus during the second semester of 2016 in a section with a total of 54 students enrolled.}
\label{ClassAttendance}
\end{figure}

Pre-class preparation of students is key to effectively profit from the in-class activities. Figure \ref{CumPreClassPreparation} shows the percentage of students that watched an insufficient, intermediate or sufficient number of videos throughout the semester as pre-class preparation. 
$50\%$ of the students did not meet the requirement of watching a sufficient number of videos, i.e. 50\% of the students studied less than 60 videos as pre-class preparation, and 30\% of the students watched an insufficient number of videos. We think that this may have had an impact on class attendance.
Figure \ref{ClassAttendance} pictures for each class the number of students that attended that class. It shows that class attendance decreased during the semester, following the same decreasing trend as for video views in Figure \ref{PreClassPreparation}. We believe, that students that were not prepared knew that the material was available online, and did not see an immediate benefit from attending class. 
They may get a wrong sense of security and confidence when having permanent access to the learning material, believing that they will find the time to acquire the knowledge by themselves. But we have observed, that generally those students have difficulties to adjust to self-regulated learning.


We did not have any mini-quiz or evaluation at the beginning of each class to check if students watched the videos (some institutions do \citep{BV}), since we wanted to enforce the fact that students should be responsible learners by their second university year. We also wanted to keep the evaluation system the closest to the traditional evaluation system for this course. 



\subsection{In-class activities and dynamics}\label{In-Class}

There are two key elements when flipping a class (see the review by \citet{Flaherty2015}), one is a successful pre-class preparation, i.e. that students study the videos consistently; and the second is to align in-class activities with pre-class activities, i.e. the instructor needs to prepare material for students to work with in class, that is in accordance with what students had prepared, that allows interaction between students and the instructor, and that permits the instructor to give valuable feedback and receive feedback on students' performance. 

We designed a worksheet for each class session, containing exercises that gradually increased their difficulty to finally achieve a general learning outcome (see \ref{..} of the supple). The difficulty of the first exercises on the worksheet was directly aligned with the video content that students had to study in their pre-class preparation. The general learning outcome as well as specific learning objectives were indicated on the worksheet for the students to know.  Students were divided in groups of 5 to 6 to solve the worksheet each class, actively collaborating and discussing among members of their own group. The participants of each group were established at the beginning of the semester and maintained through it. Table \ref{Evaluations} shows the weights of each evaluation in the course, which consisted on three tests, three quizzes and weekly worksheets.

\begin{table}[h!]
\tbl{Students in both, the flipped section and the traditional section took three common tests and three quizzes. Students in the flipped section had to turn in worksheets. Those evaluations were weighted to give the exam presentation grade (EPG). All students took a common final exam which together with the EPG gave the final grade.\\}
{\begin{tabular}{lcc} \midrule
& \multicolumn{2}{c}{Weight (\%) on exam presentation grade (EPG)}\\[.2cm]
Type of evaluation & Flipped section & Traditional section \\ \midrule
Test 1 & 25 & 25\\
Test 2 & 25 & 25\\
Test 3 & 25 & 25\\
Quiz average \textsuperscript{a} & 10 & 25\\
Worksheet average \textsuperscript{b} & 15 & 0\\ \hline \\
Final Grade & \multicolumn{2}{c}{$0.75 \times \hbox{EPG} +0.3 \times \hbox{Final exam grade}$}\\ \bottomrule
\end{tabular}}
\tabnote{\textsuperscript{a}Students took three quizzes in the semester. \textsuperscript{b} Students turned in all worksheets ones a week.}
\label{Evaluations}
\end{table}

Each group had to turn in their solved worksheet once a week, which was assigned a grade. Occasionally, students were asked to present an exercise at the blackboard and a class discussion was initiated by it. Table \ref{Activities} shows the weekly activities for students and instructors to perform during class, and their weekly duration.


\begin{table}[h!]
\centering
\tbl{ Weekly working activities of students and instructor during in-class time in lectures and help-session.\\}
{\begin{tabular}{l l  l  l  } \midrule
Type of class & Duration (min)  & Students' activity & Instructor's activity \\
& per week & & \\
\midrule	
Lecture &  3 $\times$ 70 & Peer assisted learning &  Answer doubts,  \\
& & solving a worksheet. & coach students, \\
&&&give and receive feedback.\\[.2cm]
Help-session & 70  & Students finish and turn  & Assist students in the  \\
& &in the week's worksheets. & writing of the solutions of \\
&&&exercises in the worksheets.\\ \bottomrule
\end{tabular}}
\label{Activities}
\end{table}

The effectiveness of the in-class activities gets reduced by the significant number of students that do not watch the videos before class (see 
Subsection \ref{Pre-class}). We opted by showing each class the video slide of the video that they should have watched,  as a review for students that did watch it and as a quick explanation for students that did not prepare it. This was many times not enough for the latter to be prepared to work on the designed activities, and to contribute to the team-work. The heterogeneous level of preparation slowed down the pace of the class activities and made it difficult  to guide students to achieve the learning outcome the instructor had planned.

We observed a specific behavior of the students that may depend on the fact that the topic is mathematics. We realized that students needed some time to think by themselves about the questions on the worksheet, and were not able to immediately start a discussion in their group. This produced silent moments during class that may have felt as inefficient  moments for students that did not complete their pre-class preparation, especially for less proactive students.  

As the diagram in Figure \ref{Dynamics} shows, a key consequence of in-class student centered activities is feedback. We think about feedback in two directions:  that students obtain feedback from their instructor on their progress and that the instructor obtains feedback on students' advances and pace in reaching the desired learning outcomes. 
After the individual thinking process, students realized their weaknesses in the understanding of the video topic and often needed help in order to start solving the exercises. For a numerous class this moment may be challenging if the instructor is the only available person to answer questions. It may take long for the instructor to solve all doubts and students may get demotivated and frustrated by the lack of support to advance in solving their worksheet. After the initial coaching, some groups moved forward by themselves discussing the problems and finishing the worksheet. For other groups it was difficult to move forward without any constant feedback from the instructor, especially if those groups were composed by students that lacked pre-class preparation. Usually, groups that did not finish the worksheet during class divided among themselves the exercises in order to finish the worksheet more efficiently at home. This behavior was against the purpose of the planned learning process, since the worksheet was designed such that students work gradually from the beginning to the end of the worksheet, to experience the gradual degree of difficulty of the exercises, achieving the specific learning outcomes one by one. 

The natural consequence of each class after receiving feedback is the adjustment period (see Figure \ref{Dynamics}), also in two directions:  students adjust their pre-class preparation time to be better prepared and the instructor adjusts his class to student's pace by adapting subsequent worksheets to student's weaknesses. Unfortunately, we observed that during the adjustment period often students make decisions that are counterproductive to their mastering of content, such as not attending class because of lack of preparation. We think that student adjustment actions are one of the hardest aspects to control in this process.

 We realized that in-class activities and dynamics of this type require additional support-staff with the capacity of solving students questions during in-class activities. Also, real time support to register assistance and pre-class preparation data is needed. This way the instructor could be better informed of the level of preparedness of students and their attitude towards the class format, and adapt the in-class activities to it. Therefore, to implement Flipped Learning well, there is a higher cost involved for the institution to hire more supporting-staff. This makes the flipped class format a more expensive alternative to the traditional expository class format.

\section{Results}\label{Results}

\subsection{Student performance}

\begin{figure}[htbp]
\begin{center}
\includegraphics[scale=0.65]{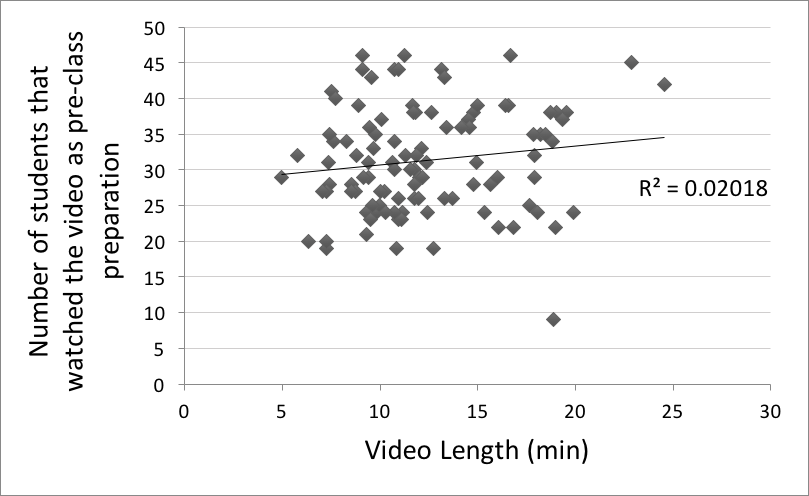}
\end{center}
\caption{The figure shows  the length of the videos in minutes vs the number of students that watched the video as pre-class preparation, with a correlation coefficient of determination $R^2=0.02018$.}
\label{LenghtWatch}
\end{figure}

 Even though the literature recommends short videos  \citep{Bergmann2012}, we did not find a clear correlation between the length of the videos and the students' interest in watching them as pre-class preparation. This  can be seen  in Figure  \ref{LenghtWatch}, where we show the video length vs the number of students that watched videos of that specific length as pre-class preparation, and a low coefficient of determination $R^2=0.02018$. If students watched the videos depended in a greater manner on the time they had available during the semester and the intersection with time consuming activities from other courses they were enrolled in.
 
\begin{figure}[htbp]
\begin{center}
\includegraphics[scale=0.6]{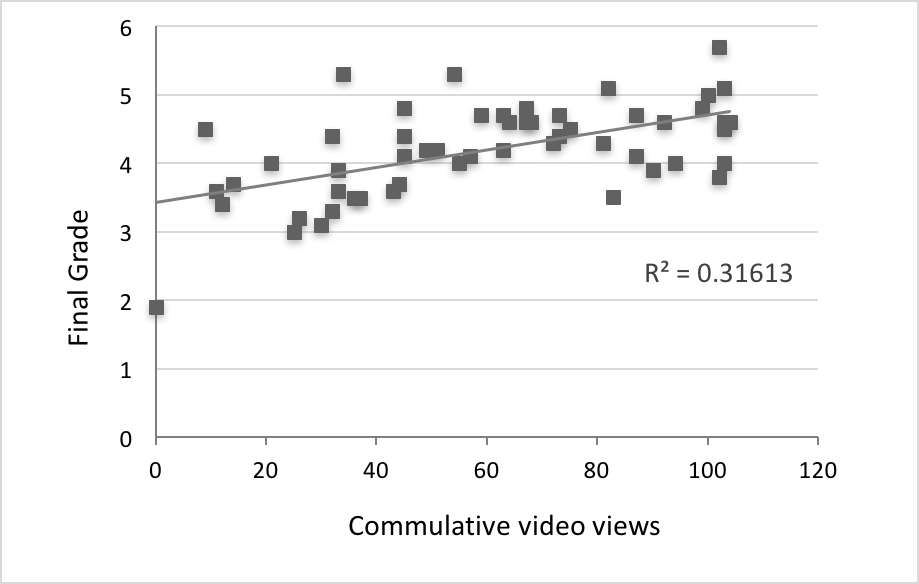}
\end{center}
\caption{The figure shows the correlation between the total number of videos a student watched during the semester as pre-class preparation (cumulative video views) vs the final grade of that student, with a coefficient of determination $R^2=0.31613$.}
\label{ViewsVsFinal}
\end{figure}

On the other hand, we did observe a correlation between the total number of videos a student watched during the semester as pre-class preparation (cumulative video views)  and that student's final grade, as can be seen in Figure \ref{ViewsVsFinal}, where we obtained a higher coefficient of determination $R^2=0.31613$.

\begin{figure}[htbp]
\begin{center}
\includegraphics[scale=0.6]{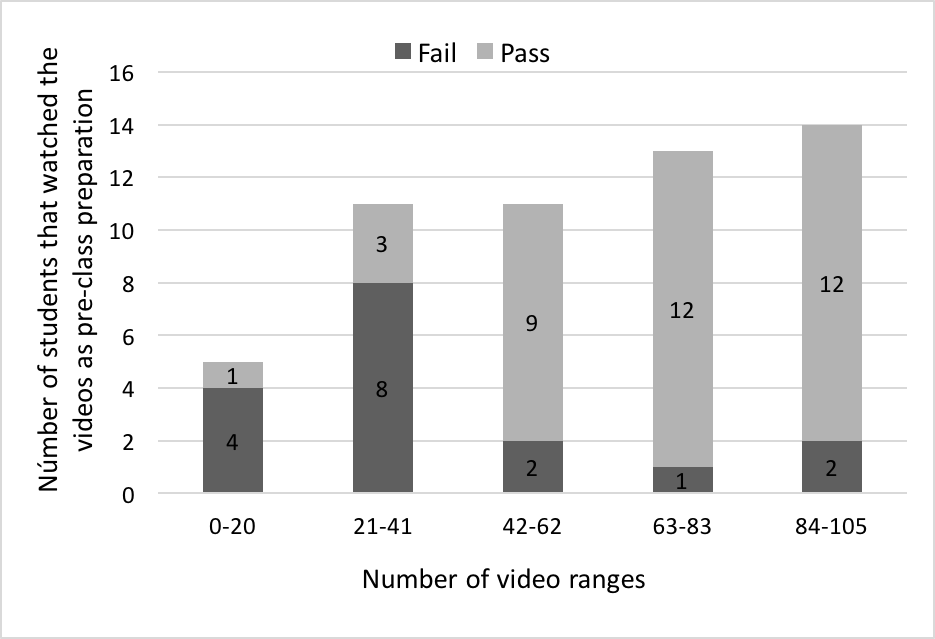}
\end{center}
\caption{Each bar represents a range of number of videos out of a total of 105 vs the number of students that watched that amount of videos in each range as pre-class preparation. Each bar shows in dark grey the number of students that failed the course and in light grey the number of students that passed the course, for each video number range setting.}
\label{Success}
\end{figure}

To give more insight we depicted in Figure \ref{Success} the performance of students in the course, measured by failing (final grade $<$ 4.0) or passing (final grade $\geq$ 4.0) the course, with respect to how many videos they had watched throughout the semester as pre-class preparation. One can observe that most of the students that watched 41 videos (39\% of the videos) or less  as pre-class preparation failed the class, but on the contrary, that most of the students that watched more than 41 videos passed the class.  Hence, most of the students that consistently prepared for class passed the course. This confirms the statement that pre-class preparation, and the compromise that comes with it,  is key for success.  

\newpage
\subsection{Flipped section vs traditional section}

\begin{figure}[htbp]
\begin{center}
\includegraphics[scale=0.6]{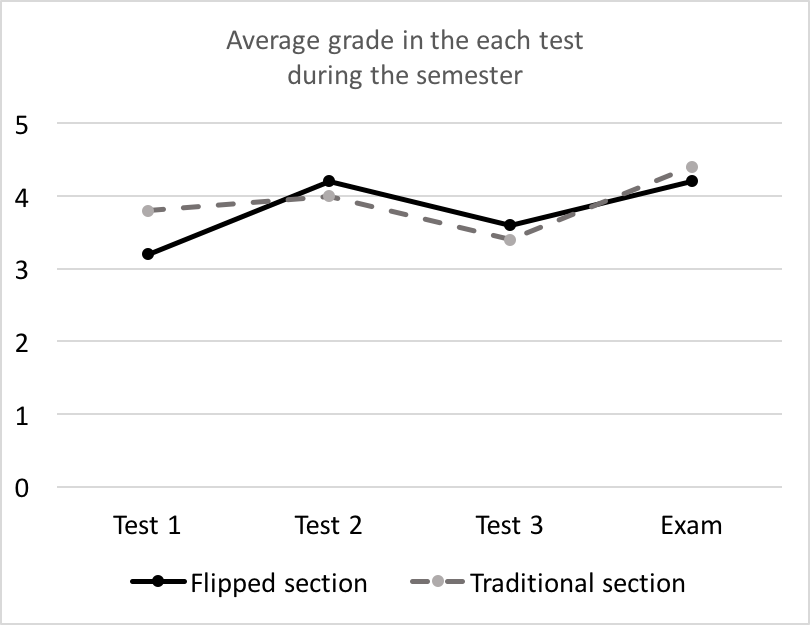}
\end{center}
\caption{The figure shows for the three common tests and final exam the average grades for the flipped section (represented with the solid line) and for the traditionally thought section (represented with the dashed line), both lectured in the second semester of 2016. The grade scale is from 1.0 to 7.0, with 4.0 being the minimum grade to pass the course.}
\label{TestGrades}
\end{figure}

We compared our results with the parallel section that thought the same course but with the traditional expository methodology, also during the second semester of 2016. Both sections had common tests and a final exam.  Figure \ref{TestGrades} shows the evolution of the average grades of those evaluations throughout the semester for both sections: The solid line represents the results for the flipped section and the dashed line the results for the traditionally thought section. 
It shows that there is not a significant difference in average test and exam results between the flipped section and the traditionally thought section. Table \ref{Indicators} shows additionally that the final grade average, standard deviation and passing/failing percentage for both sections are very similar (see Table \ref{Evaluations} for the computation of the final grade using the tests and exam grades).

\begin{table}[h!]
\centering
\tbl{ The table shows the indicators: Number of students, final grade average, standard deviation, $\%$ of students that failed and $\%$ of students that passed the course, for both, the flipped and the traditionally thought section, in the second semester of the year 2016.\\}
{\begin{tabular}{l c c}\midrule
Indicator  &  Flipped section & Traditional section\\ \midrule
Number of students & 54 & 37\\
Final grade average & 4.2 & 4.18\\
Standard deviation & 0.671 & 0.825\\
\% of students that passed & 70\% & 62\%\\
\% of students that failed & 30\% & 38\%
  \\ \bottomrule
\end{tabular}}
\label{Indicators}
\end{table}

\begin{figure}[htbp]
\begin{center}
\includegraphics[scale=0.55]{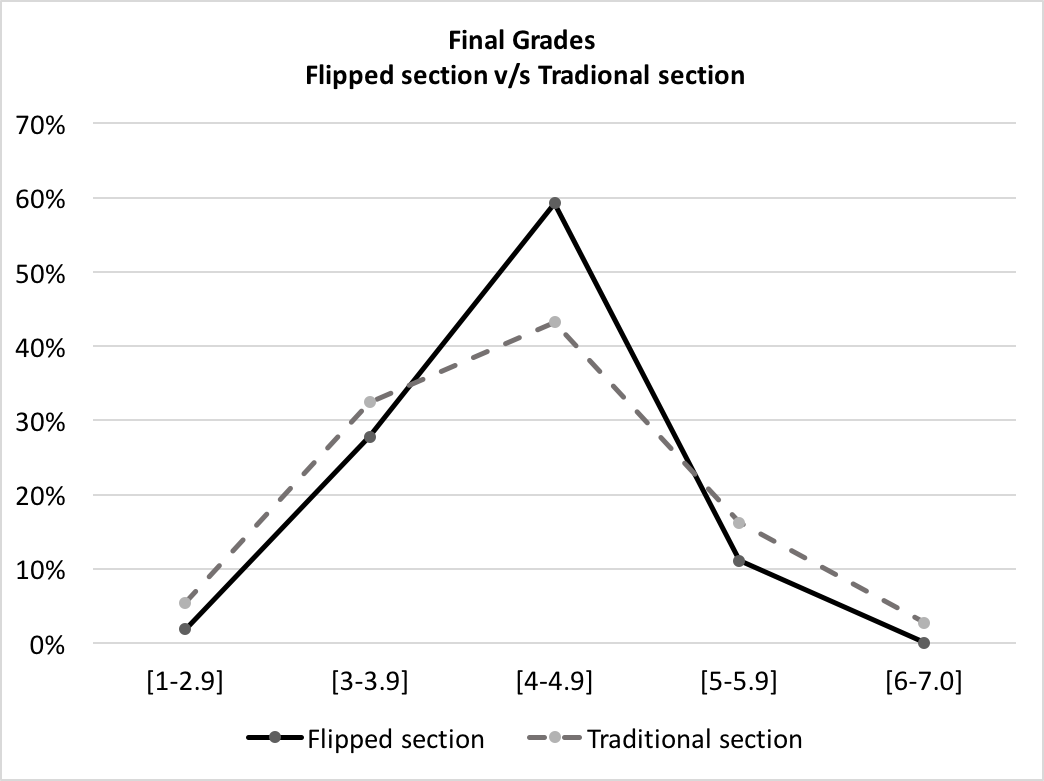}
\end{center}
\caption{The figure shows the distribution of the average final grades for the flipped section (solid line) and the traditionally thought section (dashed line). Five final grade ranges are considered: 1-2.9, 3-3.9, 4-4.9, 5-5.9 and 6-7.0.}
\label{FinalGrades}
\end{figure}

Figure \ref{FinalGrades} shows the distribution of the average final grades for the flipped section and the traditionally thought section. This gives a detailed idea of the data given in Table \ref{Indicators}. One can observe that the curves are similarly shaped following a normal distribution.

\newpage

\subsection{Student perception}

Tables \ref{PA} and \ref{NA} show qualitative results of a free text respond survey that students took at the end of the semester in the flipped section. The tables show comments that we repetitively  found in the answers and the percentages represent the percentage of students that wrote the respective comment. 

Table \ref{PA} shows positive comments. $48\%$ of the students surveyed, think that it is positive to  have access to the study material, i.e. videos and worksheets, wherever they are and whenever they need it. $48\%$ also find it positive that this methodology keeps them constantly studying. $19\%$ think that the videos are well made, and a less percentage value that this teaching form improves their learning autonomy and helps them to develop other soft skills. 
%
%


\begin{table}[h!]
\centering
\tbl{Positive aspects that appeared  repetitively in a free text respond survey taken by students at the end of the second semester 2016 in the flipped section. \\}
{\begin{tabular}{lc} \midrule
 Student's comment &   \% of students that mention this \\
  & in their free responses\\  \midrule
I can access the material when I need it & 48\%  \\
The work scheme keeps me constantly studying &  48\%\\
The videos are well made & 19\% \\
Improve my autonomy to learn & 12\%\\
Helped me develop soft skills & 10\%
  \\ \bottomrule
\end{tabular}}
\label{PA}
\end{table}

\begin{table}[h!]
\centering
\tbl{Negative aspects that appeared  repeatedly in a free text respond survey taken by students at the end of the second semester 2016 in the flipped section.}
{\begin{tabular}{lc} \midrule
 Student's comment &   \% of students that mention this \\
  & in their free responses\\  \midrule
It is a time demanding methodology  & 67\%  \\
Too many videos & 24\% \\
Videos are too long &  21\%\\
Cannot interact with the instructor while studying  & 7\%\\
Small feedback with our tasks & 7\%
  \\ \bottomrule
\end{tabular}}
\label{NA}
\end{table}




On the contrary, students also analyzed negative aspects as can be seen in Table \ref{NA}. $67\%$ of the students think that learning this way demands more of their time than the traditional expository class setting, while $24\%$ and $21\%$ of the students also think that the number of videos was too large and too long respectively. A low percentage commented that this methodology doesn't permit to interact with the instructor while studying and that not enough feedback was given.




\subsection{Methodological evolution through the years}

The methodology changed from traditional expository classes in 2015 and first semester of 2016, to Flipped Learning in the second semester of 2016 and both semesters of 2017, to finally reach a mixed methodology (between elements of Flipped Learning and traditional expository learning) in both semesters of 2018. 
In this section we will provide some results on the evolution of the following three aspects in time between the years 2015 and 2018: First, student's choice to enroll in a flipped class section; second, instructor's evaluation in the course by semester; and third, the average percentage of students that pass the course. 

\begin{figure}[htbp]
\begin{center}
\includegraphics[scale=0.5]{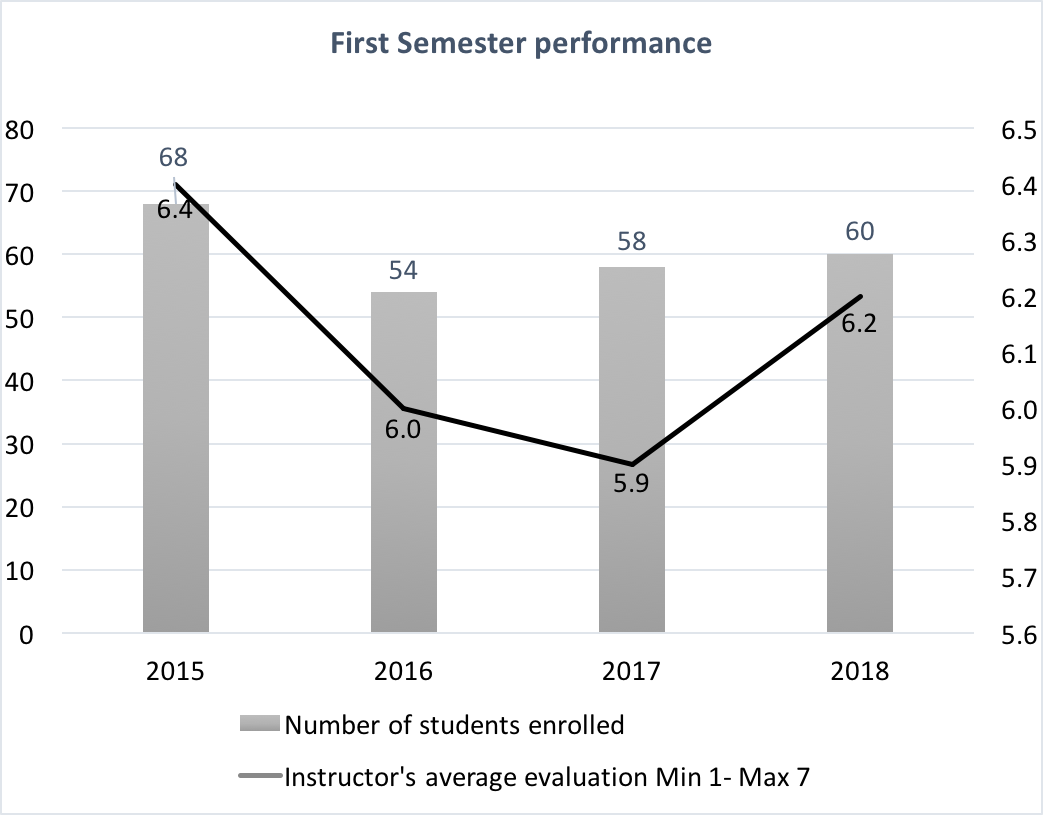}
\end{center}
\caption{The figure shows for the first semester of each of the years 2015 to 2018, the number of students enrolled (bars with respect to the left vertical axis) and the instructor's evaluation (solid line with respect to the right vertical axis) in a scale from 1.0 to 7.0, being 1.0 the worst and 7.0 be best. }
\label{EnrollEvalS1}
\end{figure}

\begin{figure}[htbp]
\begin{center}
\includegraphics[scale=0.5]{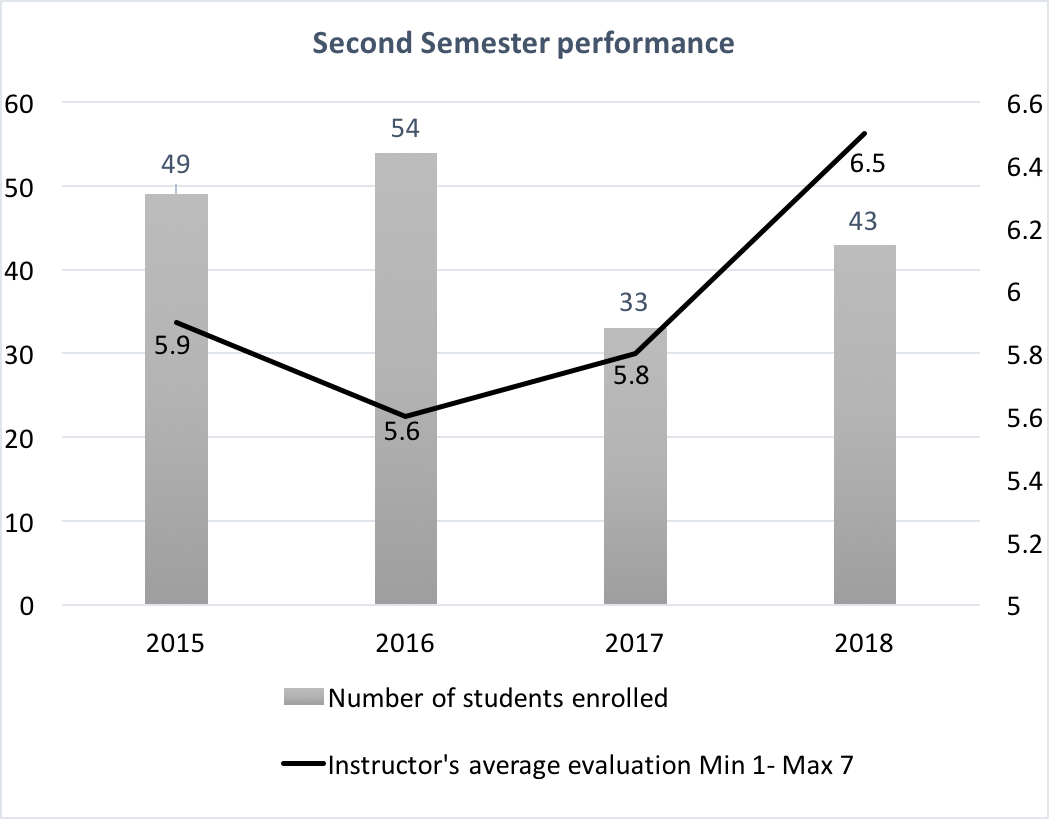}
\end{center}
\caption{The figure shows for the second semester of each of the years 2015 to 2018, the number of students enrolled (bars with respect to the left vertical axis) and the instructor's evaluation (solid line with respect to the right vertical axis) in a scale from 1.0 to 7.0, being 1.0 the worst and 7.0 be best. }
\label{EnrollEvalS2}
\end{figure}

The bars in Figure \ref{EnrollEvalS1} and \ref{EnrollEvalS2} show the number of students that enrolled in the course, taught by the same professor, the first semester and the second semester respectively of each academic year  between 2015 and 2018.  

 Students that are behind in their academic progress and failed some mathematics course before enrolling in Multivariable Calculus, generally take this course in the second semester (hence would correspond to the enrollment shown in Figure  \ref{EnrollEvalS2}), and are usually weaker students academically. Figure \ref{EnrollEvalS2} shows a low enrollment (33 students) in the second semester of 2017 compared to other semesters. We attribute this to the information that different generations of students shared with each other after the implementation of Flipped Learning the two previous semesters (second semester of 2016 and first semester of 2017), informing students about to enroll that flipped courses require more work.

The line in Figure \ref{EnrollEvalS1} and \ref{EnrollEvalS2} shows the average instructor's evaluation each semester, for the same professor, for the first semester and second semester respectively, of each academic year between 2015 and 2018. The evaluation range is from 1.0 to 7.0, being 1.0 the worst and 7.0 the best. In both figures one can see that students give the professor a significantly lower evaluation in semesters where flipped methodology was implemented (second semester of 2016 (Figure \ref{EnrollEvalS2}), and first \& second semester of 2017 (Figure \ref{EnrollEvalS1} and \ref{EnrollEvalS2})), compared to the other semesters where the course was thought using the traditional expository methodology or the mixed methodology. Also, especially students that are repeating the course (second semester enrollment, Figure \ref{EnrollEvalS2}) evaluate the professor the best when mixed methodology is used (see Figure \ref{EnrollEvalS2} year 2018), and the same students evaluate the professor the worst in the semester that were thought using Flipped Learning (see Figure \ref{EnrollEvalS2} years 2016, 2017).

\begin{figure}[htbp]
\begin{center}
\includegraphics[scale=0.47]{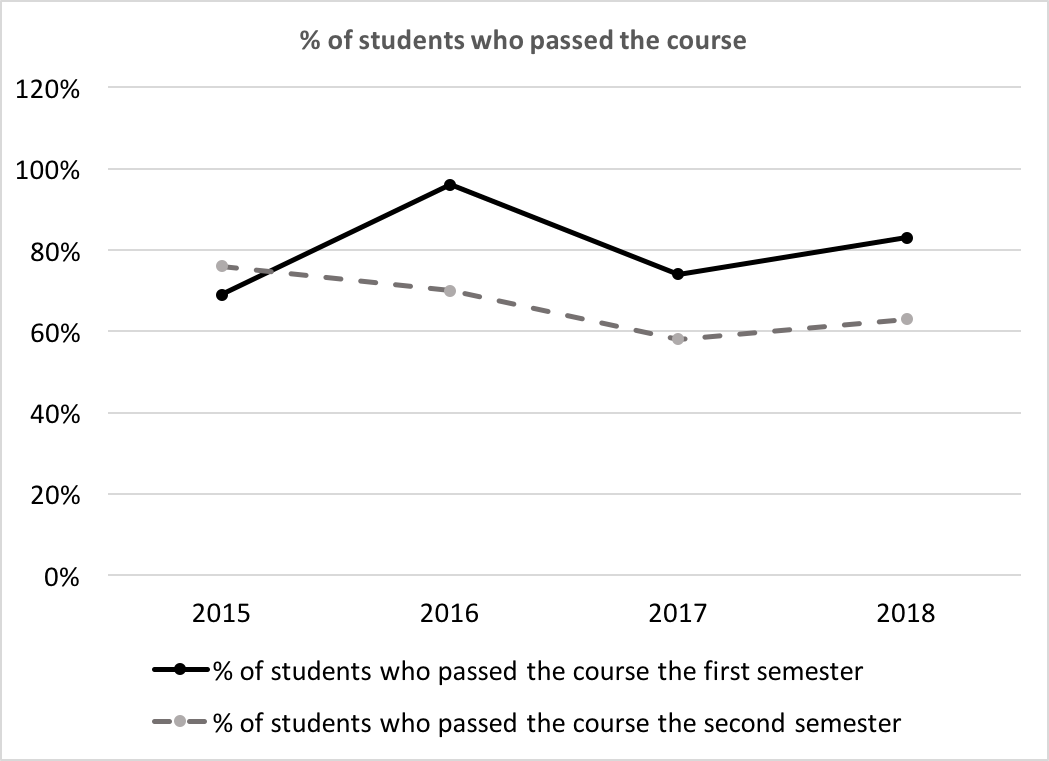}
\end{center}
\caption{The figure shows the percentage of students that passed the course throughout the years 2015 to 2018. The solid line represents first semesters and the dashed line second semesters results.}
\label{Passing}
\end{figure}

Finally, Figure \ref{Passing} shows the percentage of students that passed the course throughout the years 2015 to 2018, where the solid line represents first semesters and the dashed line second semesters. One can see that the passing percentage is generally higher during first semesters, which is when mainly non repeating students are enrolled, than second semesters which is when the course is generally attended by repeating students. Also, for both, repeating students and non-repeating students, the passing percentage decreases when the flipped methodology is used (second semester 2016 and first \& second semester 2017). 

The professor that thought this course throughout the years 2015 to 2018, also observed a consistent decrease in attendance throughout all semesters in which Flipped Learning was implemented, similar to what was observed during the first implementation in the second semester of 2016 (see Figure \ref{ClassAttendance}). In general, we could observe that if students have the choice to enroll in a traditional classroom or a flipped classroom they generally will choose the first.

\newpage

\section{Discussion and conclusion: Flipped Learning may not be the answer but it helps to move towards a modern way of teaching}\label{DisConcl}

After our experience of creating and implementing an active methodology, we asked ourselves: Why is it justified to invest so much effort and time in designing and implementing Flipped Learning? In this section we want to discuss answers to this question.


Despite experiencing some difficulties in implementing Flipped Learning, such as students not watching the videos before in-class activities, or difficulties during face-to-face in class dynamics, this experience taught us a way of teaching that adjusts better to the present time, where technology and online communication is present everywhere. 
We found that Flipped Learning creates a space to develop a closer student-teacher and student-student relationship that allows free discussions and interactions to clarify questions and fortify weaknesses, as well as improving soft skills. Unfortunately, that space cannot be used in the most efficient way, reaching all students, if certain conditions are not satisfied. These conditions may seem obvious, but are necessary to keep in mind when planing a teaching innovation involving Flipped Learning. The first condition is that one needs to take into account that the pace of a flipped class is often determined by how fast students reach the desired learning outcome. Hence, if students posses a weak academic  background, it may be necessary to reduce the number of topics the class covers. The Multivariable Calculus course we flipped, covers topics such as Taylor's Theorem, Surface Integrals, Divergence and Stokes Theorems. These topics are not easy to understand, therefore the course thought in an active learning section may easily move at a much slower pace than a parallel section where the topics are transmitted in an expository way, with less student-instructor interaction. A second condition for successfully using the space Flipped Learning creates, is to count with a support staff. For instance, a teaching assistant that records students' assistance and video views before class, and who also helps during in-class activities with the coaching of students is essential. A third condition is that in order  to keep track of students progress and provide valuable feedback to each individual, it is crucial to have small class sizes. Otherwise it becomes difficult to evaluate the advances of the class, class pace and contents, to reach the desired learning outcomes.

%

 
 
 Our results did not show significant differences in student's performance with respect to the parallel section taught in a traditional expository manner, nevertheless we think that the ideal methodology for our students and the resources available is not to completely flip a class or to teach the expository traditional way, but to combine methodologies (mixed methodology), especially for students that are repeating the course. Up to which measure Flipped Learning activities will work, depends on the characteristics of the group of students taking the class, with more proactive group of students reacting better to the flipped classroom than others. This shows that instructors are not replaceable, since some students need elements of the traditional expository class setting and guiding in order to be successful learners.


During our trajectory as instructors, we implemented several methodological changes pushed by the idea of making our classes more appealing and to obtain better learning outcomes. Usually, to implement those changes required much energy and resources which led to not keeping those projects running. This time, our experience was different since even though the project still required much time and energy, it left us with a final product that we can use in the future in multiple ways, depending on the methodology we want to implement. We are now implementing a mixed methodology, that combines flipped-learning elements with traditional expository classes, depending on the topics covered and how  receptive the group of students is to active learning. We post all videos onto a platform available to students, as well as worksheets containing exercises to achieve mastery in the topics treated in each video. Each student can watch the videos according to their own time and pace, either as preparation for class or as a class review if the topic was covered during an expository lecture. Students appreciate and take good advantage of the video resources in their own way. If they do not watch the videos as preparation, they feel that class-time will not be lost since expository elements are always included.  Even if they miss a class, they feel included in the progress of the course because of the video material that they can watch to keep track of the topics covered. In general, topics that are difficult for students or too long to cover in detail during class time are put in a video capsule, giving the lecturer extra time during class to exercise base on those difficult sometimes abstract mathematics topics. This way, students as well as professors have an extra tool to achieve students' mastery in the course.

Finally, we also discovered that this experience was an important step towards creating online courses. We feel better prepared to develop an online learning platform for our institution, or participating in the course offer at platforms such as Coursera \citep{Co}. We now know the difficulties in creating audiovisual material and the requirements to succeed in finishing such a project successfully. We understand our personal weaknesses and institutional weaknesses that need to be addressed when moving into methodologies using digital elements, which we think is inevitable in our time. It is interesting and challenging to notice, that even though the way of thinking and developing basic mathematics hasn't changed, we feel the need to teach it in a very different way compared to how they did  hundred of years ago. The way students learn is changing.

\section*{Acknowledgement(s)}
We thank colleagues from the Engineering and Sciences Faculty University Adolfo Ib\'a\~nez (UAI) for their valuable feedback and especially our collegues in the Mathematics Department UAI. We thank members of the Educational Learning Center UAI for fruitful discussions.

\section*{Funding}
This research was partially funded by``The Teaching Innovation Competition 2016" of the Educational Learning Center and the Engineering and Science Faculty of the Adolfo Ib\'a\~nez University in Santiago de Chile.


\end{document}